\documentclass{gtmon_a}
\pdfoutput=1

\usepackage{pinlabel}

\proceedingstitle{Heegaard splittings of 3--manifolds (Haifa 2005)}
\conferencestart{10 July 2005}
\conferenceend{19 July 2005}
\conferencename{Heegaard splittings of 3--manifolds}
\conferencelocation{Haifa}

\editor{Cameron Gordon}
\givenname{Cameron}
\surname{Gordon}

\editor{Yoav}
\givenname{Yoav}
\surname{Moriah}

\title{Destabilizing amalgamated Heegaard splittings}

\author{Jennifer Schultens} 
\givenname{Jennifer}
\surname{Schultens}
\address{Department of Mathematics\\
University of California, Davis\\\newline
1 Shields Avenue\\
Davis CA 95616\\
USA}
\email{jcs@math.ucdavis.edu}

\author{Richard Weidmann} 
\givenname{Richard}
\surname{Weidmann}
\address{Department of Mathematics and Maxwell Institute of Mathematical
Sciences\\\newline
Heriot-Watt University\\
Riccarton\\
Edinburgh\\
EH14 4AS\\
Scotland}
\email{R.Weidmann@ma.hw.ac.uk}

\volumenumber{12}
\issuenumber{}
\publicationyear{2007}
\papernumber{13}
\startpage{319}
\endpage{334}

\doi{}
\MR{}
\Zbl{}

\arxivreference{math.GT/0510386}


\received{20 January 2007}
\revised{12 March 2007}
\accepted{12 March 2007}
\published{3 December 2007}
\publishedonline{3 December 2007}
\proposed{}
\seconded{}
\corresponding{}
\version{}

\makeatletter
\def\cnewtheorem#1[#2]#3{\newtheorem{#1}{#3}
\expandafter\let\csname c@#1\endcsname\c@thm}

\let\xysavmatrix\xymatrix
\def\xymatrix{\disablesubscriptcorrection\xysavmatrix}
\AtBeginDocument{\let\bar\wbar\let\tilde\wtilde\let\hat\what}

\theoremstyle{plain}
\newtheorem{thm}{Theorem}
\cnewtheorem{lem}[thm]{Lemma}
\cnewtheorem{cor}[thm]{Corollary}

\theoremstyle{definition}
\cnewtheorem{dfn}[thm]{Definition}
\cnewtheorem{defn}[thm]{Definition}
\cnewtheorem{rem}[thm]{Remark}
\cnewtheorem{notation}[thm]{Notation}
\cnewtheorem{convention}[thm]{Convention}
\cnewtheorem{exa}[thm]{Example}
\cnewtheorem{exercise}[thm]{Exercise}

\makeatother  %  move after \newtheorem block

\makeautorefname{defn}{Definition}
\makeop{genus}

\begin{document}

\begin{webabstract}
We construct a sequence of pairs of 3--manifolds $(M_1^n, M_2^n)$
each with incompressible torus boundary and with the following two properties: 

(1) For $M^n$ the result of a carefully chosen glueing of $M_1^n$ and
$M_2^n$ along their boundary tori, the genera $(g_1^n, g_2^n)$ of
$(M_1^n, M_2^n)$ and the genus $g^n$ of $M^n$ satisfy the inequality

$g^n/(g_1^n + g_2^n) < 1/2.$

(2) The result of amalgamating certain unstabilized Heegaard
splittings of $M_1^n$ and $M_2^n$ to form a Heegaard splitting of
$M$ produces a stabilized Heegaard splitting that can be
destabilized successively $n$ times.
\end{webabstract}

\begin{htmlabstract}
<p class="noindent">
We construct a sequence of pairs of 3&ndash;manifolds (M<sub>1</sub><sup>n</sup>, M<sub>2</sub><sup>n</sup>)
each with incompressible torus boundary and with the following two properties:
</p>
<p class="noindent">
(1) For M<sup>n</sup> the result of a carefully chosen glueing of M<sub>1</sub><sup>n</sup> and
M<sub>2</sub><sup>n</sup> along their boundary tori, the genera (g<sub>1</sub><sup>n</sup>, g<sub>2</sub><sup>n</sup>) of
(M<sub>1</sub><sup>n</sup>, M<sub>2</sub><sup>n</sup>) and the genus g<sup>n</sup> of M<sup>n</sup> satisfy the inequality
<br>
g<sup>n</sup>/(g<sub>1</sub><sup>n</sup> + g<sub>2</sub><sup>n</sup>) < 1/2.
</p>
<p class="noindent">
(2) The result of amalgamating certain unstabilized Heegaard
splittings of M<sub>1</sub><sup>n</sup> and M<sub>2</sub><sup>n</sup> to form a Heegaard splitting of
M produces a stabilized Heegaard splitting that can be
destabilized successively n times.
</p>
\end{htmlabstract}

\begin{abstract}
We construct a sequence of pairs of 3--manifolds $(M_1^n, M_2^n)$
each with incompressible torus boundary and with the following two properties: 

(1)\qua For $M^n$ the result of a carefully chosen glueing of $M_1^n$ and
$M_2^n$ along their boundary tori, the genera $(g_1^n, g_2^n)$ of
$(M_1^n, M_2^n)$ and the genus $g^n$ of $M^n$ satisfy the inequality
\[\frac{g^n}{g_1^n + g_2^n} < \frac{1}{2}.\]
(2)\qua The result of amalgamating certain unstabilized Heegaard
splittings of $M_1^n$ and $M_2^n$ to form a Heegaard splitting of
$M$ produces a stabilized Heegaard splitting that can be
destabilized successively $n$ times.
\end{abstract}

\begin{asciiabstract}
We construct a sequence of pairs of 3-manifolds (M_1^n, M_2^n) each
with incompressible torus boundary and with the following two
properties:

(1) For M^n the result of a carefully chosen glueing of M_1^n and
M_2^n along their boundary tori, the genera (g_1^n, g_2^n) of
(M_1^n, M_2^n) and the genus g^n of M^n satisfy the inequality
g^n/(g_1^n + g_2^n) < 1/2.

(2) The result of amalgamating certain unstabilized Heegaard
splittings of M_1^n and M_2^n to form a Heegaard splitting of
M produces a stabilized Heegaard splitting that can be
destabilized successively n times.
\end{asciiabstract}

\maketitle

\section{Introduction}

About 10 years ago, Cameron McA Gordon asked the following question:
Can the pairwise connect sum of two 3--manifolds each with an
unstabilized Heegaard splitting yield a 3--manifold with a stabilized
Heegaard splitting?  This question stumped the experts for many years
 but recently a negative answer to this question
has been announced by D Bachman \cite{B} and R Qiu \cite{Q}.

More generally, one can ask how Heegaard splittings behave under other
types of ``sums'', that is, when the 3--manifolds containing them are
glued along positive genus boundary components.  How Heegaard genus
behaves under these circumstances is one of the many questions
investigated by Klaus Johannson in \cite{Jo} and by the first author
in \cite{Sch}.  In both cases, inequalities relating the Heegaard
genus of the glued 3--manifold to the Heegaard genera of the original
3--manifolds are obtained.  Most strikingly, the inequalities give
lower bounds on the Heegaard genus of the glued 3--manifold in terms of
the Heegaard genera of the original 3--manifolds.  But these lower
bounds are fractions of the sum of the genera of the original
3--manifolds.  A better bound under more restrictive circumstances has
recently been obtained by D Bachman, E Sedgwick and S Schleimer
\cite{BSS}.

One upshot is that, in general, the phenomenon of ``degeneration of
Heegaard genus'' under glueing of 3--manifolds can't be ruled out.  It
appears however that under certain, possibly generic circumstances, this
phenomenon does not occur.  For instance, in \cite{La} Marc
Lackenby shows that for a pair of hyperbolic 3--manifolds each with one
boundary component and under certain restrictions on the glueing,
minimal genus Heegaard splittings of the glued 3--manifold are always obtained from
Heegaard splittings of the original 3--manifolds by amalgamation.

It is presently unknown how large ``degeneration of Heegaard genus''
under glueing can be.  Interestingly, the issue of stabilization
implicitly arises in the investigation of this phenomenon in \cite{SS}
and in \cite{Sch}.  The examples given in this note make this issue explicit.  In
particular, we provide examples that illustrate how ``degeneration of
Heegaard genus'' under glueing corresponds to
the existence of stabilizations in the amalgamation of Heegaard
splittings of the original 3--manifolds.  In doing so, we provide
counterexamples to a conjecture of Kobayashi, Qiu, Rieck and Wang \cite{KQRW}.

\subsection*{Acknowledgments}
The first author is partially supported by the
grant DMS-0603736 from the National Science Foundation.

\section{Definitions}

For standard definitions and results concerning knots, see Burde and
Zieschang \cite{BZ}, Lickorish \cite{L} or Rolfsen \cite{R}.  For standard
definitions and results pertaining to $3$--manifolds, see Hempel \cite{H}
or Jaco \cite{J}.

\begin{defn}
A {\em height function} on ${\mathbb S}^3$ is a Morse function with exactly
two critical points.  
\end{defn}

This last assumption guarantees that $h$ induces
a foliation of $S^3$ by spheres, along with one maximum that we denote
by $\infty$ and one minimum that we denote by $-\infty$.

\begin{defn}
Let $K$ be a knot in $S^3$. If all minima of $h|_K$ occur below all maxima of $h|_K$, then we say
that $K$ is in {\em bridge position} with respect to $h$.  The {\em bridge number}
of K, $b(K)$, is the minimal number of maxima required for $h|_K$.
\end{defn}

\begin{defn}
If $K$ is in bridge position, then a regular level surface below
all maxima and above all minima is called a {\em bridge
surface}.  
\end{defn}

\begin{defn}
An {\em upper disk (lower disk)} is an embedded disk whose interior is disjoint from the knot whose boundary is partitioned
into two subarcs, one contained in a bridge surface and one a subarc
of the knot that lies above (below) the bridge surface.  A {\em strict
upper disk (strict lower disk)} is an upper (lower) disk whose interior
lies above (below) the bridge surface.

A {\em complete set of strict upper (lower) disks} is a set of upper (lower) disks
such that each subarc of the knot lying above (below) the bridge
surface meets exactly one disk in the set.
\end{defn}

\begin{defn} \label{defn:cb} 
A {\em compression body} is a $3$--manifold $W$ obtained from $S\times I$ where $S$ is a
closed orientable connected surface by attaching $2$--handles to $S \times
\{0\} \subset S \times I$ and capping off any resulting $2$--sphere
boundary components with $3$--handles.  We denote $S \times \{1\}$ by
$\partial_+W$ and $\partial W - \partial_+W$ by $\partial_-W$.
Dually, a compression body is an orientable $3$--manifold obtained from
a closed orientable surface $\partial_-W \times I$ or a 3--ball or a
union of the two by attaching $1$--handles.

In the case where $\partial_-W = \emptyset$, we also call $W$ a
handlebody.  
\end{defn}

\begin{defn}
Let ${\cal A}=\{a_1,\ldots ,a_k\}$ be a collection of annuli in a compression body
$W$.  Then ${\cal A}$ is a {\em primitive collection} if there is a
collection ${\cal D}=\{D_1,\ldots ,D_k\}$ of pairwise disjoint properly embedded  disks in $W$ such that
$a_i$ meets $D_i$
in a single spanning arc and $a_i\cap D_j=\emptyset$ for $j\neq i$.
\end{defn}

\begin{defn} A set of {\em defining disks} for a compression body
$W$ is a set of disks $\{D_1, \dots, D_n\}$ properly embedded in $W$
with $\partial D_i \subset \partial_+W$ for $i = 1$, $\dots, n$ such
that the result of cutting $W$ along $D_1 \cup \dots \cup D_n$ is
homeomorphic to $\partial_-W \times I$ or a $3$--ball in the case that
$W$ is a handlebody.
\end{defn}

\begin{defn} \label{defn:Heegaard splitting} A {\em Heegaard
splitting} of a $3$--manifold $M$ is a pair $(V, W)$ in which $V$, $W$
are compression bodies and such that $M = V \cup W$ and $V \cap W =
\partial_+V = \partial_+W = S$.  We call $S$ the {\em splitting surface} or
{\em Heegaard surface}.  Two Heegaard splittings are considered {\em equivalent}
if their splitting surfaces are isotopic.

The {\em genus} of $M$, denoted by $g(M)$, is the smallest possible genus of
the splitting surface of a Heegaard splitting for $M$.  
\end{defn}

\begin{defn}
Let $(V, W)$ be a Heegaard splitting.  A Heegaard splitting is
called {\em stabilized} if there is a pair of properly embedded disks $(D, E)$ with $D \subset V$ and
$E \subset W$ such that $\#|\partial D \cap \partial E | =
1$.  We call the pair of disks $(D, E)$ a {\em stabilizing pair of disks}.  A
Heegaard splitting is {\em unstabilized} if it is not stabilized.
\end{defn}

\begin{defn}
{\em Destabilizing} a Heegaard splitting $(V, W)$ is the act of creating a
Heegaard splitting from $(V, W)$ by performing ambient 2--surgery on
$S$ along the cocore of a $1$--handle in either $V$ or $W$.
\end{defn}

Note that the result of performing ambient 2--surgery on $S$ along the
cocore of a $1$--handle in either $V$ or $W$ is not necessarily a
Heegaard splitting.  In order for this operation to be a
destabilization, the result is required to be a Heegaard splitting.
This definition is equivalent to presupposing a stabilizing pair of
disks $(D, E)$ and cutting along $D$. (Here $D$ is the cocore of a
$1$--handle of $V$ and the existence of $E$ guarantees that the result
of cutting along $D$ is a Heegaard splitting.)

\begin{figure}[ht!]
\labellist\small
\pinlabel {$a_1$} [b] at 220 160
\pinlabel {$a_2$} [l] at 158 139
\endlabellist
\centerline{\includegraphics[width=0.5\textwidth]{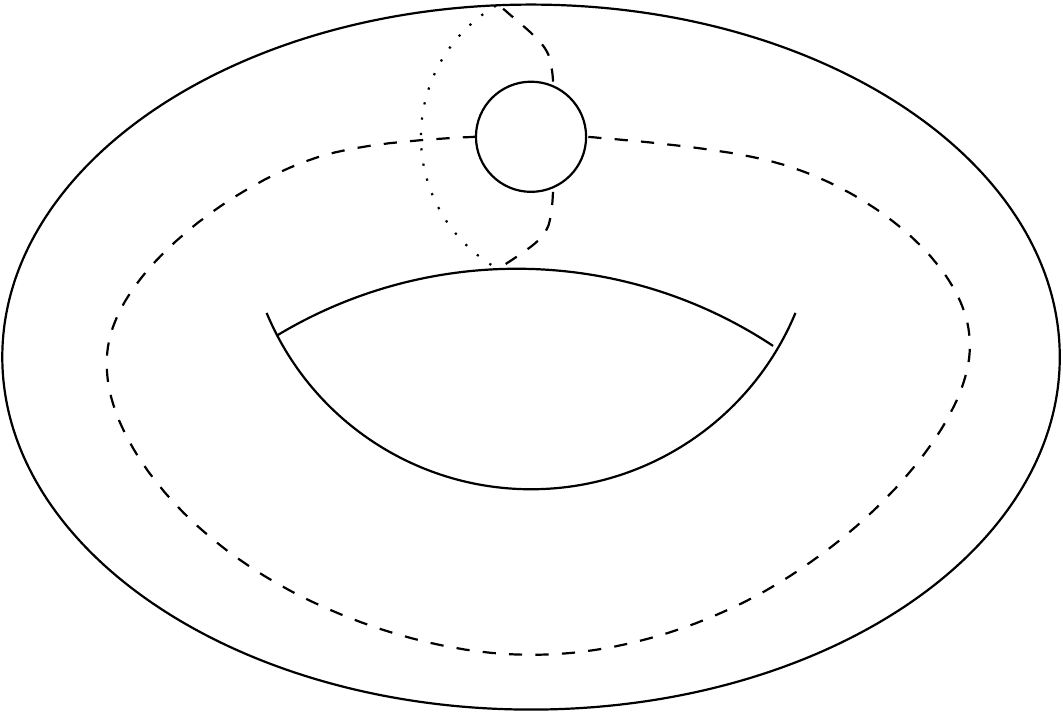}}
\caption{Arcs $a_1, a_2$ for a punctured torus}
\label{dsa1}
\end{figure}

\begin{figure}[ht!]
\labellist\small
\pinlabel {$b_1$} [t] at 40 165
\pinlabel {$b_2$} [r] at 135 245
\pinlabel {$b_3$} [t] at 235 170
\endlabellist
\centerline{\includegraphics[width=0.5\textwidth]{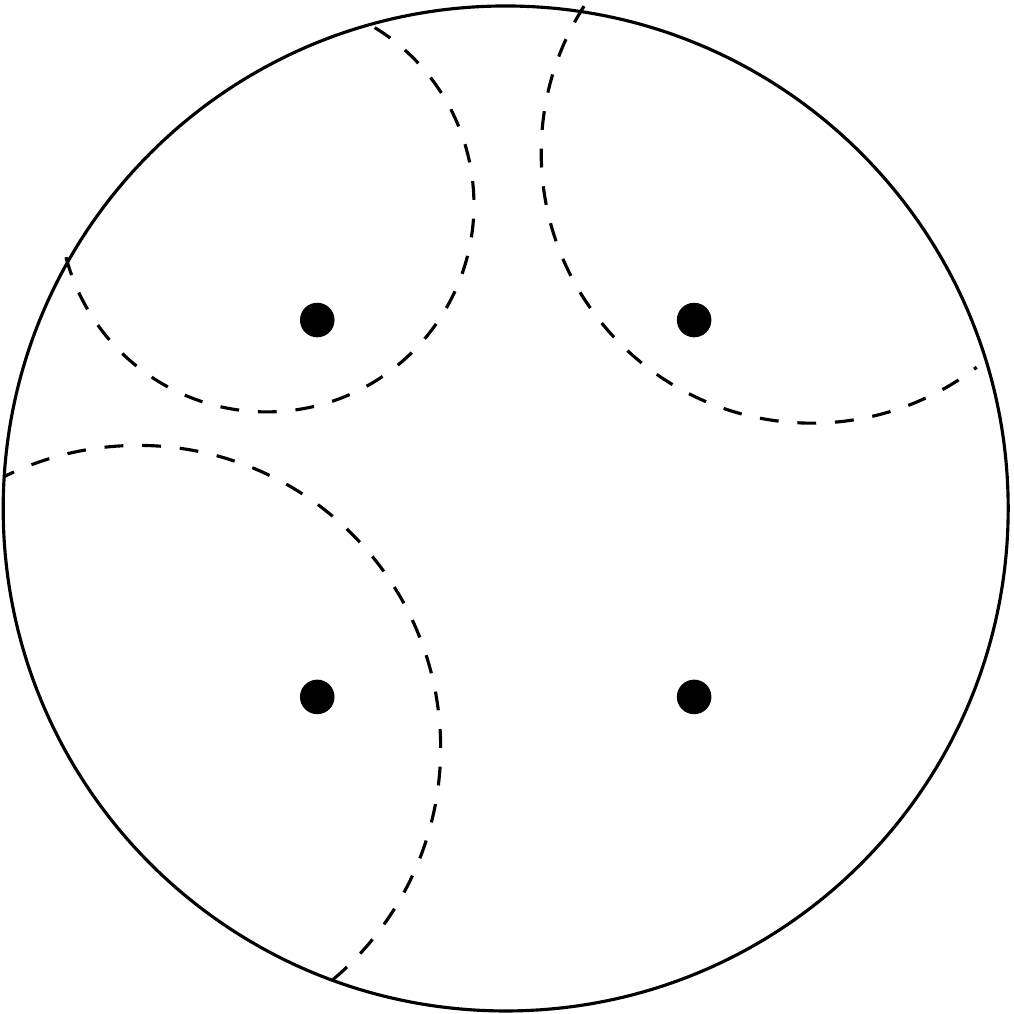}}
\caption{Arcs $b_1, b_2, b_3$ for an orbifold with four cone points}
\label{dsa7}
\end{figure}

\begin{defn} \label{stHS}
  Let $M$ be a compact orientable Seifert fibered space with quotient
  space an orientable orbifold $Q$.  Denote the genus of the surface
  underlying $Q$ by $g$ and the number of cone points by $n$.
  Assume further that $M$ (and hence $Q$) has exactly one boundary
  component.  (This simplifying assumption is met in all examples
  considered here.)
  
  Let $a_1, \dots, a_{2g}, b_1, \dots, b_{n-1}$ be a disjoint
  collection of arcs in $Q$ that cut $Q$ into disks each containing at
  most one cone point.  In the case of the once punctured torus, such
  a collection of arcs is shown in \fullref{dsa1}.  In the case of
  an orbifold with underlying surface a disk and with four cone
  points, such a collection of arcs is shown in \fullref{dsa7}.
  If the underlying surface of $Q$ is a disk, we further assume that
  each arc $b_i$ cuts off a subdisk from $Q$ containing exactly one
  cone point.  

Abusing notation slightly, denote a collection of arcs in $M$ that
  projects to $a_1, \dots, $ $a_{2g}, b_1, \dots, b_{n-1}$ also by
  $a_1, \dots, a_{2g}, b_1, \dots, b_{n-1}$.  Now take $V$ to be a
  regular neighborhood of $a_1, \dots, a_{2g}, b_1, \dots, b_{n-1}$
  together with a regular neighborhood of $\partial Q \times {\mathbb
  S}^1$.  Take $W$ to be the closure of the complement of $V$ in $M$.
  It is an easy exercise to show that $(V, W)$ is a Heegaard splitting
  of $M$.  Such a Heegaard splitting is called a {\em vertical} Heegaard
  splitting of $M$. If $Q$ has no cone points, that is, if $M=Q\times
  S^1$, then this splitting is also called the {\em standard} Heegaard
  splitting of $M$.
  \end{defn}

\begin{defn}
A {\em tunnel system} for a knot $K$ in ${\mathbb S}^3$ is a collection of
arcs $t_1, \dots, t_n$ such that the complement of $K \cup t_1 \cup
\dots \cup t_n$ is a handlebody.  The {\em tunnel number} of a knot $K$ is
the least number of arcs required for a tunnel system of $K$.
\end{defn}

\begin{defn}
Suppose $K$ is in bridge position and that there are $n$ maxima.  We
may assume temporarily that all maxima occur in the same level surface
$L$.  The maxima may be connected by a system of $n - 1$ disjoint arcs
in $L$.  It is an easy exercise to show that this set of arcs is a
tunnel system.  It is called an {\em upper tunnel system}.

The same exercise shows that there is a set of defining disks ${\cal
D}$ for the complement of $K \cup t_1 \cup \dots \cup t_n$ of the
following type: Each component of ${\cal D}$ has interior below $L$,
furthermore, below $L$, its boundary runs once along exactly one
component of $K - K \cap L$.  This set of disks is called a {\em complete
set of lower disks} for the upper tunnel system. 
\end{defn}

\begin{defn}
Suppose $t_1, \dots, t_n$ is a tunnel system for a knot $K$ in
${\mathbb S}^3$.  Denote the complement of $K$ by $M$.  Take $V$ to be
a regular neighborhood of $\partial M \cup t_1 \cup \dots \cup t_n$
and take $W$ to the closure of the complement of $V$.  Then $(V, W)$
is a Heegaard splitting called the {\em Heegaard splitting corresponding to
the tunnel system} $t_1, \dots, t_n$.
\end{defn}

The definition of amalgamation is a lengthy one.  In the last 15
years, this term has been used in the following context: A pair of
3--manifolds $M_1, M_2$ each with a Heegaard splitting are identified
along components of their boundary.  This results in a 3--manifold $M$.
The Heegaard splittings of $M_1, M_2$ can be used to construct a
canonical Heegaard splitting of $M$ called the amalgamation of the two
Heegaard splittings.  One assumes that in each of $M_1, M_2$ the
boundary components along which the glueing occurs are contained in a
single compression body.  Roughly speaking, then, the collars of the
boundary components lying in this compression body are discarded and
the remnants of the two compression bodies in $M_1 - \text{(collars)}$
identified to the remnants of the two compression bodies in $M_2 -
\text{(collars)}$.  This is done in such a way that the 1--handles that are
attached to the collar on such a boundary component in $M_1$ become
attached to the compression body in $M_2$ that does not meet any of
the boundary components along which the glueing takes place and vice
versa.  For a formal definition see below.

\begin{defn}
Let $M_1, M_2$ be $3$--manifolds with $R$ a closed subsurface of
$\partial M_1$, and $S$ a closed subsurface of $\partial M_2$.
Suppose that $R$ is homeomorphic to $S$ via a homeomorphism $h$.
Further, let $(X_1, Y_1), (X_2, Y_2)$ be Heegaard splittings of $M_1,
M_2$.  Suppose further that $N(R) \subset X_1, N(S) \subset X_2$.
Then, for some $R' \subset \partial M_1 \backslash R$ and $S' \subset
\partial M_2 \backslash S$, $X_1 = N(R \cup R') \cup (\text{1--handles})$ and
$X_2 = N(S \cup S') \cup (\text{1--handles})$.  Here $N(R)$ is homeomorphic to
$R \times I$ via a homeomorphism $f$ and $N(S)$ is homeomorphic to $S
\times I$ via a homeomorphism $g$.  Let $\sim$ be the equivalence
relation on $M_1 \cup M_2$ generated by

(1) $x \sim y$ if $x, y \; \epsilon \; \eta(R)$ and $p_1 \cdot f(x) = p_1
    \cdot f(y)$,

(2) $x \sim y$ if $x, y \; \epsilon \; \eta(S)$ and $p_1 \cdot g(x) = p_1
    \cdot g(y)$,

(3) $x \sim y$ if $x \; \epsilon \; R$, $y \; \epsilon \; S$ and $h(x) = y$,

where $p_1$ is projection onto the first coordinate.  Perform
isotopies so that for $D$ an attaching disk for a $1$--handle in $X_1,
D'$ an attaching disk for a $1$--handle in $X_2$, $[D] \cap [D'] =
\emptyset$.  Set $M = (M_1 \cup M_2)/\sim, X = (X_1 \cup Y_2)/\sim,$
and $Y = (Y_1 \cup X_2)/\sim$.  In particular, $(N(R) \cup N(S)/\sim)
\cong R, S$.  Then $X = Y_2 \cup N(R') \cup (\text{1--handles})$, where the
$1$--handles are attached to $\partial_+Y_2$ and connect $\partial
N(R')$ to $\partial_+Y_2$.  Hence $X$ is a compression body.
Analogously, $Y$ is a compression body.  So $(X, Y)$ is a Heegaard
splitting of $M$.  The splitting $(X, Y)$ is called the {\em amalgamation}
of $(X_1, Y_1)$ and $(X_2, Y_2)$ along {R, S} via $h$.
\end{defn}

\section{A single destabilization}

We first consider a concrete example that illustrates the issues under
discussion.  Let $T_i$ be a punctured torus for $i=1,2$. As
3--manifolds $M_1, M_2$ we take $T_i \times {\mathbb S}^1$ for $i = 1,
2$.  Note that $\partial M_i$ is a torus, for $i = 1, 2$.  We take $M$
to be the result of glueing $M_1$ to $M_2$ in such a way that
$(\partial T_1) \times \{1\}$ and $(\partial T_2) \times \{p\}$ have
intersection number one on the resulting torus.

We describe two distinct Hegaard splittings for $M$:

\begin{exa} \label{g4}
Let ${\mathbb S}^1 = I_1 \cup I_2$ be a decomposition of ${\mathbb
S}^1$ into two intervals that meet at their endpoints.  Let $V_i =
T_i \times I_1$ and $W_i = T_i
\times I_2$, for $i = 1, 2$.  Then $V_i$ and $W_i$ are genus 2
handlebodies.  Denote the annulus in which $V_i$ meets
$\partial M_i$ by $A_i$ and that in which $W_i$ meets $\partial M_i$
by $B_i$.  Due to the choice of glueing of $\partial M_1$ and
$\partial M_2$ that results in $M$, $A_1$ meets $A_2$ in a (square) disk.  As
do $B_1$ and $B_2$.  In other words, $V = V_1 \cup V_2$ is
homeomorphic to the result of taking the disjoint union of $V_1$ and
$V_2$ and joining the two components by a 1--handle.  In particular, it
is a genus 4 handlebody.  The same is true for $W = W_1 \cup W_2$.
Thus $(V, W)$ is a genus 4 Heegaard splitting of $M$.
\end{exa}

\begin{exa} \label{amalg}
Let $(X_i, Y_i)$ be the standard Heegaard splitting of $M_i$ (see \fullref{stHS}), for $i =
1, 2$.  And let $(X, Y)$ be the amalgamation of $(X_1, Y_1)$ and
$(X_2, Y_2)$
\end{exa}

\begin{thm}
The genus of $M_i$ is three for $i = 1, 2$ and the genus of $M$ is
four.
\end{thm}

\begin{proof}
Recall that the rank, that is, the smallest number of generators, of the
fundamental group of a 3--manifold provides a lower bound for the genus
of a Heegaard splitting of that 3--manifold.  Here
\[\pi_1(M_i) = F_2\oplus \mathbb Z \] 
Abelianizing yields a free abelian group of rank 3.  Thus $\hbox{rank
}\pi_1(M_i)=3$ and hence the Heegaard splitting constructed in \fullref{amalg} has minimal genus.

The Seifert--Van Kampen Theorem yields a presentation of $\pi_1(M)$ as
$$\pi_1(M_1)*_{\mathbb Z^2}\pi_1(M_2).$$ Quotienting out the normal
closure of the amalgamated subgroup yields $\mathbb Z^2*\mathbb
Z^2$ as this kills the fibre and the commutator of the generating pair of the free base group on both sides. It follows that $$\hbox{rank }\pi_1(M)\ge \hbox{rank }\mathbb
Z^2*\mathbb Z^2=4.$$Hence the Heegaard splitting in \fullref{g4}
has minimal genus.
\end{proof}

The fact that the minimal genus Heegaard splitting is less than the genus of a minimal genus amalgamation in these examples illustrates a phenomenon known as ``degeneration of Heegaard genus" under glueing.

\begin{thm} \label{1destab}
The Heegaard splitting $(X, Y)$ of $M$ is stabilized.
\end{thm}

\begin{figure}%[here]
\centerline{\includegraphics[width=0.3\textwidth]{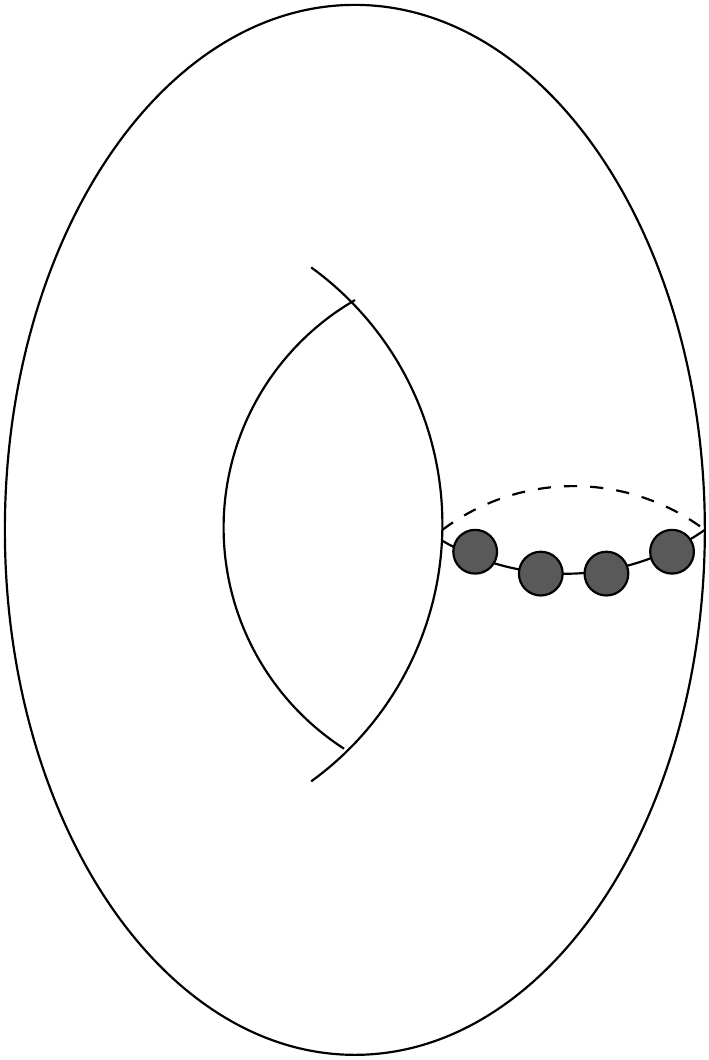}}
\caption{The boundary of $D_i$ as it appears on $\partial M_i$}
\label{dsa2}
\end{figure}

\begin{figure}%[here]
\centerline{\includegraphics[width=0.3\textwidth]{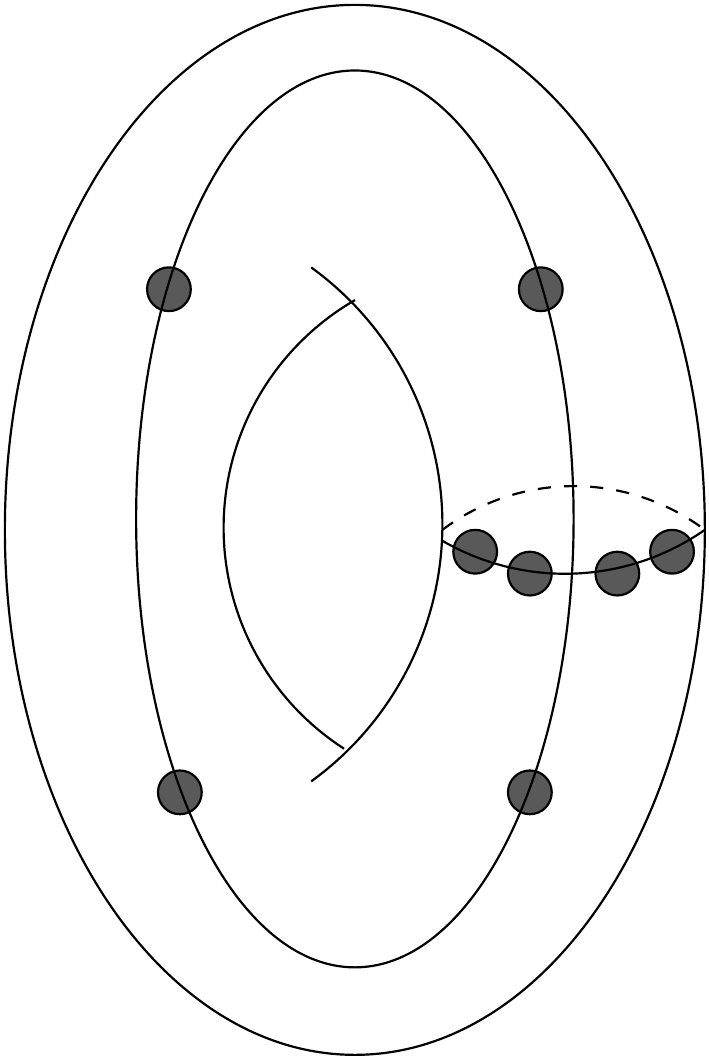}}
\caption{The boundaries of $D_1$ and $D_2$ as they intersect}
\label{dsa3}
\end{figure}

\begin{proof}
For $i=1,2$, choose arcs $a^i_1, a^2_1$ in $T_i\subset M_i$ as
in \fullref{stHS}.  Then $T_i - (N(a^i_1) \cup N(a^i_2))$ is a
disk $D_i$.  It's boundary meets $\partial M_i$ as in \fullref{dsa2}.  After the amalgamation, a copy of $D_i$ survives in $M_i
\subset M$, for $i = 1, 2$.  How $\partial D_1$ and $\partial D_2$
intersect is pictured in \fullref{dsa3}.  Thus $(D_1, D_2)$ are a
stabilizing pair of disks.
\end{proof}

\begin{cor}
The Heegaard splitting $(X, Y)$ of $M$ can be destabilized exactly
once.
\end{cor}

\begin{exercise}
Show that destabilizing the Heegaard splitting in \fullref{amalg} yields
the Heegaard splitting in \fullref{g4}.
\end{exercise}

\section{Arbitrarily many  destabilizations}

We now construct a sequence of pairs of 3--manifolds that exhibit a more
general phenomenon.  More specifically, for each $n$, we construct a pair
$(M_1^n, M_2^n)$ of 3--manifolds as follows:  Given $n$, take $M_1^n$
to be a Seifert fibered space with base orbifold a disk with $n + 1$
cone points.  We denote the natural quotient map from $M_1^n$ to the base
orbifold by $p_n$.  Take $K^n$ to be a knot that has bridge number $n$
and tunnel number $n - 1$.  The existence of such knots is guaranteed
by \cite[Theorem 0.1]{LM}.  Indeed, in \cite{LM}, M Lustig and Y Moriah
define the class of generalized Montesinos knots.  The referenced theorem
provides very technical but nevertheless achievable sufficient conditions
under which such a knot has bridge number $n$ and tunnel number $n - 1$.
Take $M_2^n$ to be the complement of $K^n$ in ${\mathbb S}^3$.

Glue $M_1^n$ to $M_2^n$ in such a way that a fiber of $M_1^n$ is identified with
a meridian of $M_2^n$.  Denote the 3--manifold obtained in this way by
$M^n$.  Consider the following Heegaard splittings of $M^n$:

\begin{exa} \label{g4n}
  Let $b_1, \dots, b_n$ be a collection of arcs that cut the base
  orbifold of $M_1^n$ into disks each with exactly one cone point.
  Bicolor these disks red and blue, that is, color these disks in such
  a way that disks abutting along an arc are given distinct colors.
  The preimage of these arcs in $M_1^n$ is a collection of annuli that
  cut $M_1^n$ into solid tori.  These tori inherit colors from the
  bicoloring of the disks to which they project.  Take $V_1^n$ to be
  the union of the red tori and $W_1^n$ to be the union of the blue
  tori.

Let $P$ be a bridge sphere for $K^n$.  Then $P$ divides $M_2^n$ into
two components that we label $V_2^n$ and $W_2^n$. We can clearly
assume that the $2n$ meridional boundary curves of $P\cap M_2^n$ match
up with the boundary curves of the annuli $b_1,\ldots ,b_n$ Now set
$V^n = V_1^n \cup V_2^n$ and $W^n = W_1^n \cup W_2^n$.
\end{exa}

\begin{lem} \label{g4nhs}
The decomposition $(V^n, W^n)$ is a Heegaard splitting of $M^n$.
\end{lem}

We first prove an auxilliary lemma.  It is well known, but we include
it here for completeness.

\begin{lem} \label{glueprim}
Suppose $X$ and $Y$ are handlebodies.  Let ${\it A}$ be a collection
of $k$ essential annuli in $\partial X$ and let ${\it B}$ be a
primitive collection of $k$ annuli in $\partial Y$.  Glue $X$ to $Y$
by identifying ${\it A}$ and ${\it B}$.  Denote the result by $E$.
Then $E$ is a handlebody.
\end{lem}

\begin{proof}
Since ${\it B}$ is a primitive collection of $k$ annuli in
$\partial Y$, there is a collection ${\cal Y}$ of $k$ disjoint
essential disks such that each annulus meets one of the disks in
exactly one arc and is disjoint from the other disks.  Cutting $Y$
along ${\cal Y}$ yields a handlebody $Y'$ and cuts each component of
${\it B}$ into a disk.  The remnants of ${\cal Y} \cup {\it B}$ on
$\partial Y'$ are disks.  Thus a set of defining disks for $Y'$ can be
isotoped to be disjoint from the remnants of ${\cal Y} \cup {\it B}$
on $\partial Y'$.  Hence they can be used to augment ${\cal Y}$ to a
set of defining disks ${\cal Y}'$ of $Y$.

Choose a set of defining disks ${\cal X}$ for $X$.  We may assume that
each component of ${\cal X}$ meets each component of ${\it A}$ in
spanning arcs.  (Note that each component of ${\it A}$ is met by a non
zero number of such arcs, because it is essential.)  In $E$ we can
place a copy of the appropriate element of ${\cal Y}$ along each such
spanning arc.  Thus in $E$, the components of ${\cal X}$ can be
extended into $Y \subset E$ by parallel copies of elements of ${\cal
Y}$ to an embedded disk.  Denote the set of disks resulting from
${\cal X}$ via these extensions along with a set of defining disks for
$Y'$ by ${\cal E}$.

Now the result of cutting $E$ along ${\cal E}$ is a 3--ball since it
can also be obtained by glueing 3--balls (the result of cutting $Y$
along ${\cal E} \cap Y$) to a 3--ball (the result of cutting $X$ along
${\cal X}$) along disks.  It follows that $E$ is a handlebody.
\end{proof}

We now prove \fullref{g4nhs}.  Fortunately, the hard work has
already been accomplished.

\begin{proof}[Proof of \fullref{g4nhs}]
To see that $(V^n, W^n)$ is a Heegaard
splitting, consider the following: Each component of $V_1^n$ and
$W_1^n$ is a solid torus.  In particular, it is a handlebody.
Furthermore, both $V_2^n$ and $W_2^n$ are genus $n$ handlebodies each
meeting $\partial M_2^n$ in a primitive collection of $n$ annuli.  More
specifically, we can take a complete set of strict upper disks or a
complete set of strict lower disks, respectively, to be the required
collection of disks.  See \fullref{dsa8}.

\begin{figure}[ht!]
\centerline{\includegraphics[width=0.7\textwidth]{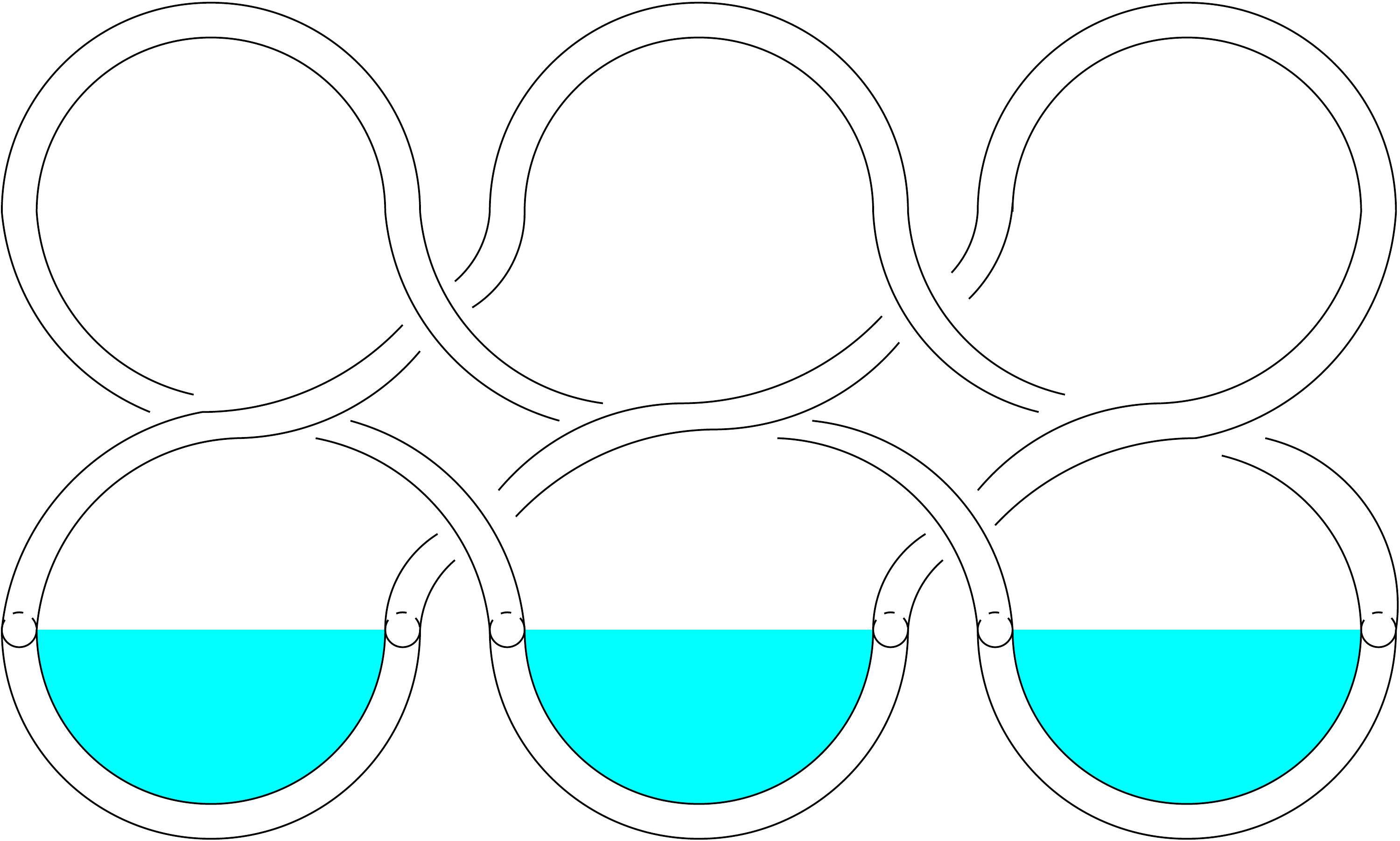}}
\caption{The submanifold $V_2^1$ or $W_2^1$ of $M_2^1$ with
a collection of disks meeting primitive annuli as required}
\label{dsa8}
\end{figure}

It thus follows
from \fullref{glueprim} that $V^n$ and $W^n$ are handlebodies.  Thus
$(V^n, W^n)$ is a Heegaard splitting.
\end{proof}

\begin{exa} \label{amalgn}
Take $(X_1^n, Y_1^n)$ to be a vertical Heegaard splitting of
$M_1^n$.  Take\break $t_1, \dots, t_{n-1}$ to be an upper tunnel system of
$M_2^n$ and take $(X_2^n, Y_2^n)$ to be the Heegaard splitting
corresponding to $t_1, \dots, t_{n-1}$.  Now take $(X^n, Y^n)$ to be
the Heegaard splitting of $M^n$ resulting from the amalgamation of
$(X_1^n, Y_1^n)$ and $(Y_1^n, Y_2^n)$.
\end{exa}

\begin{thm}
For $M_1^n, M_2^n, M^n$ as above,
\[\genus(M_1^n) + \genus(M_2^n) - \genus(M^n) > n\]
and
\[\frac{\genus(M^n)}{\genus(M_1^n) + \genus(M_2^n)} < \frac{1}{2}\]
\end{thm}

\begin{proof}
Note first that $\pi_1(M_1^n)$ maps onto the fundamental group of the
base orbifold which is a free product of $n+1$ cyclic groups. Thus
$\pi_1(M_1^n)$ is of rank $n+1$ by Grushko's theorem. It follows that
the genus of $M_1^n$ is $n+1$.  Furthermore, since the tunnel number
of $K^n$ is $n-1$, the genus of $M_2^n$ is $n$.  The Heegaard splitting
constructed in \fullref{g4n} bears witness to the fact that the
Heegaard genus of $M^n$ is at most $n$.
\end{proof}

Again, the manifold pairs $M_1^n, M_2^n$ exhibit the phenomenon of
``degeneration of Heegaard genus'' under glueing.

Note that the genus of a Heegaard splitting of $M^n$ resulting from an
amalgamation of minimal genus Heegaard splittings is $2n$.  In particular,
the genus of $(X^n, Y^n)$ is $2n$.

\begin{thm}
There are $n$ disjoint pairs of stabilizing disks for $(X^n, Y^n)$.  In
other words, the Heegaard splitting $(X^n, Y^n)$ of $M^n$ can be
destabilized successively at least $n$ times.  Specifically, the
Heegaard splitting obtained from $(X^n, Y^n)$ is the result of
stabilizing $(V^n, W^n)$ $n$ times.
\end{thm}

\begin{proof}
Recall that $M_2^n$ is the complement of $K^n$ and that $Y_2^n$ is the
complement of $K^n$ together with an upper tunnel system.  See \fullref{dsa9}.  Recall also that after amalgamation, (a collar of)
$Y_2^n$ is a subset of $X^n$.

\begin{figure}[ht!]
\centerline{\includegraphics[width=0.7\textwidth]{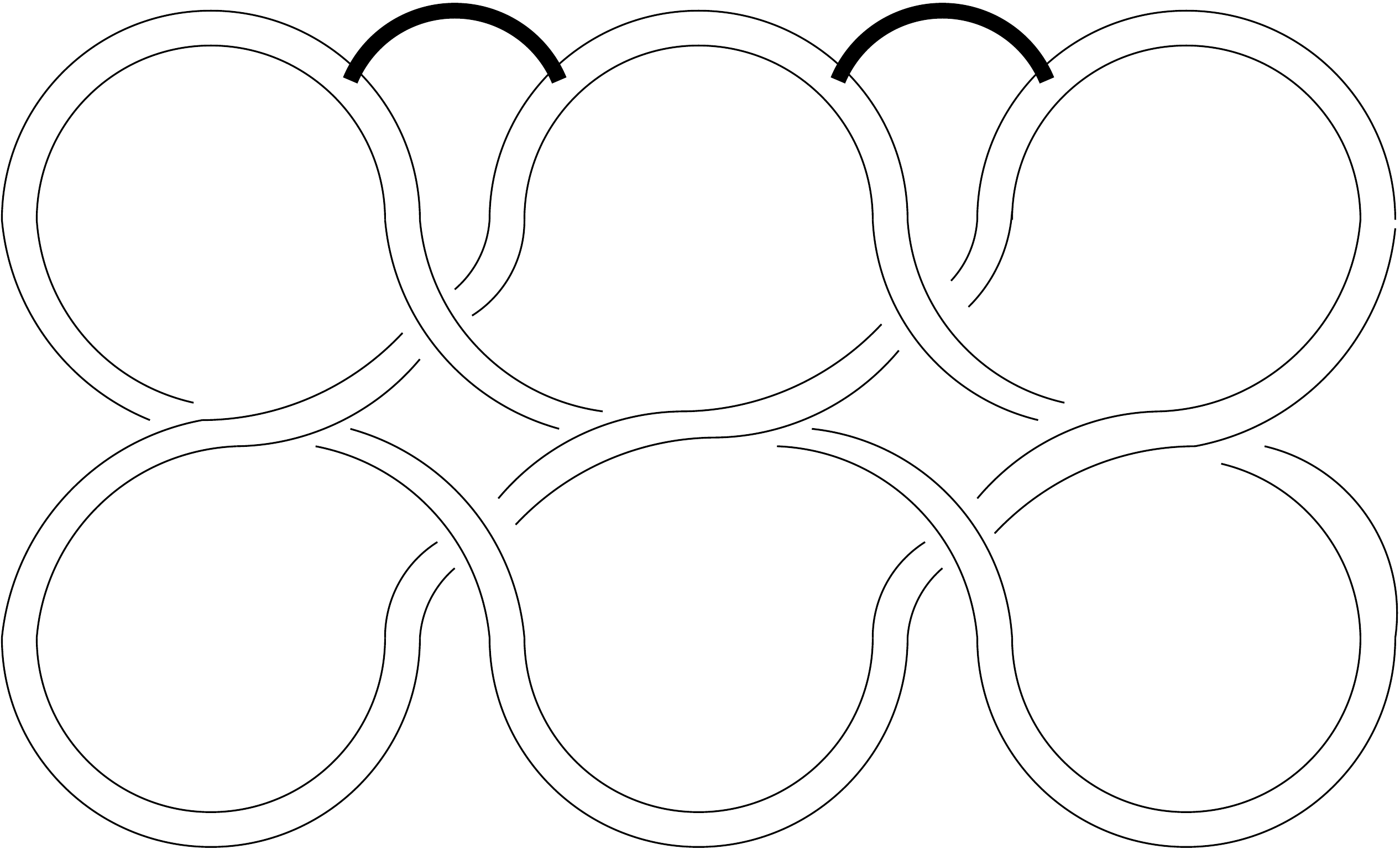}}
\caption{$K^3$ with an upper tunnel system}
\label{dsa9}
\end{figure}

Denote the torus resulting from the identification of $\partial M_1^n$
and $\partial M_2^n$ by $T$.  Recall that after the amalgamation, the
torus $T$ minus the attaching disks for the 1--handles with cores $b_1,
\dots, b_n$ to one side and the upper tunnel system to the other
side lies in the splitting surface $F^n$ of $(X^n, Y^n)$. We isotope
$n$ essential subannuli of $T$ into $M_1^n$ and denote the resulting
annuli by $U_1, \dots, U_n$.  We isotope the other $n$ subannuli of
$T$ into $M_2^n$ and denote the result by $A_1, \dots, A_n$.  We
subdivide $T$ into these subannuli in such a way that $U_1, \dots,
U_n$ are vertical in $M_1^n$ and $A_1, \dots, A_n$ are meridional in
$M_2^n$.  Furthermore, we subdivide $T$ into these subannuli in such a
way that $U_i$ meets the endpoints of exactly two distinct components
of $b_1, \dots, b_n$.  See Figures~\ref{dsa10} and \ref{dsa11}.

\begin{figure}[ht!]
\labellist\small
\pinlabel {$b_1$} [t] at 55 165
\pinlabel {$b_2$} [r] at 148 237
\pinlabel {$b_3$} [t] at 250 170
\endlabellist
\centerline{\includegraphics[width=0.5\textwidth]{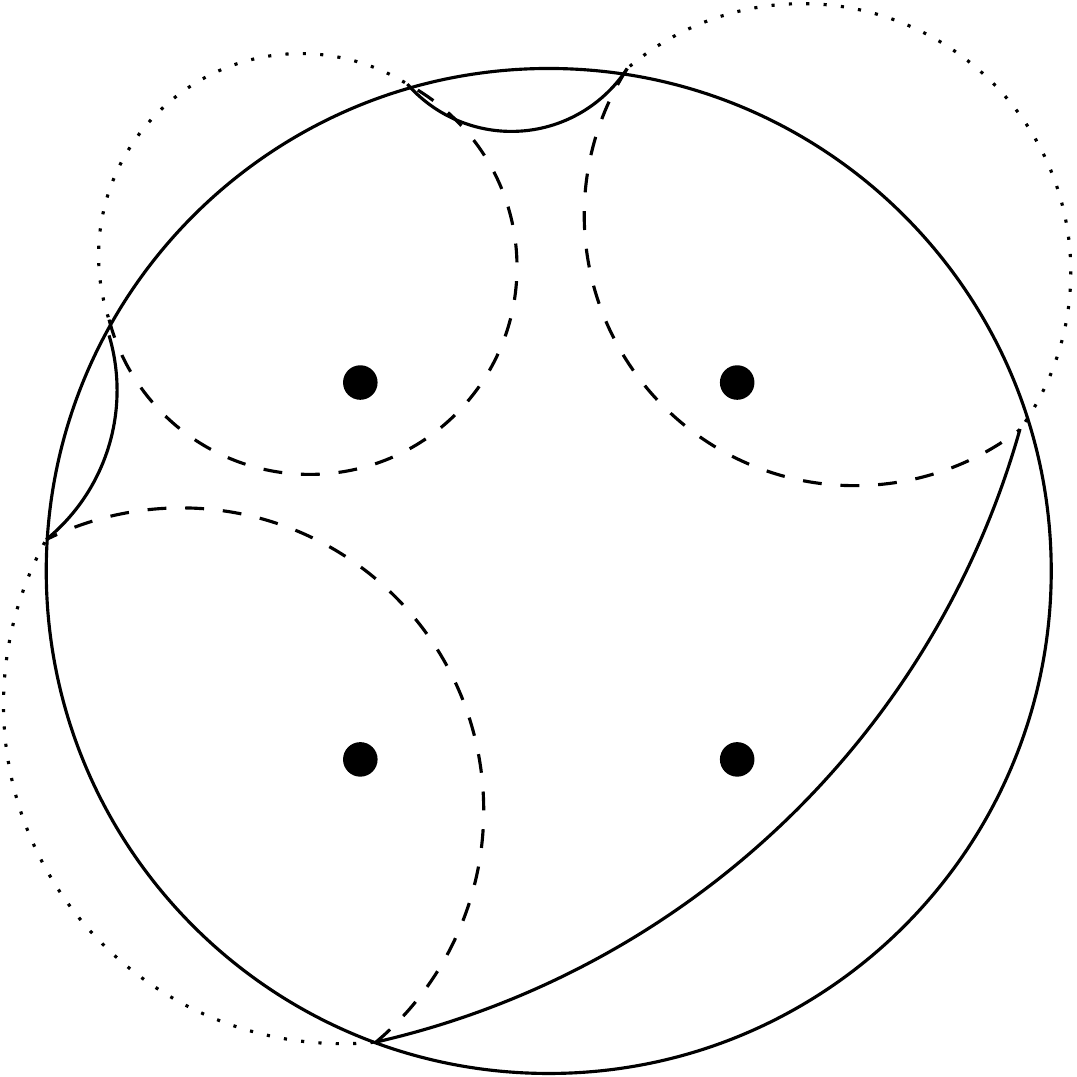}}
\caption{Vertical annuli in $M_1^3$}
\label{dsa10}
\end{figure}

\begin{figure}[ht!]
\labellist\small
\pinlabel {tunnel tube} [b] at 233 496
\pinlabel {tunnel tube} [b] at 505 496
\pinlabel {$A_1$} [b] at 108 480
\pinlabel {$A_2$} [b] at 365 488
\pinlabel {$A_3$} [b] at 613 474
\pinlabel {$T$} [b] at 360 20
\endlabellist
\centerline{\includegraphics[width=0.7\textwidth]{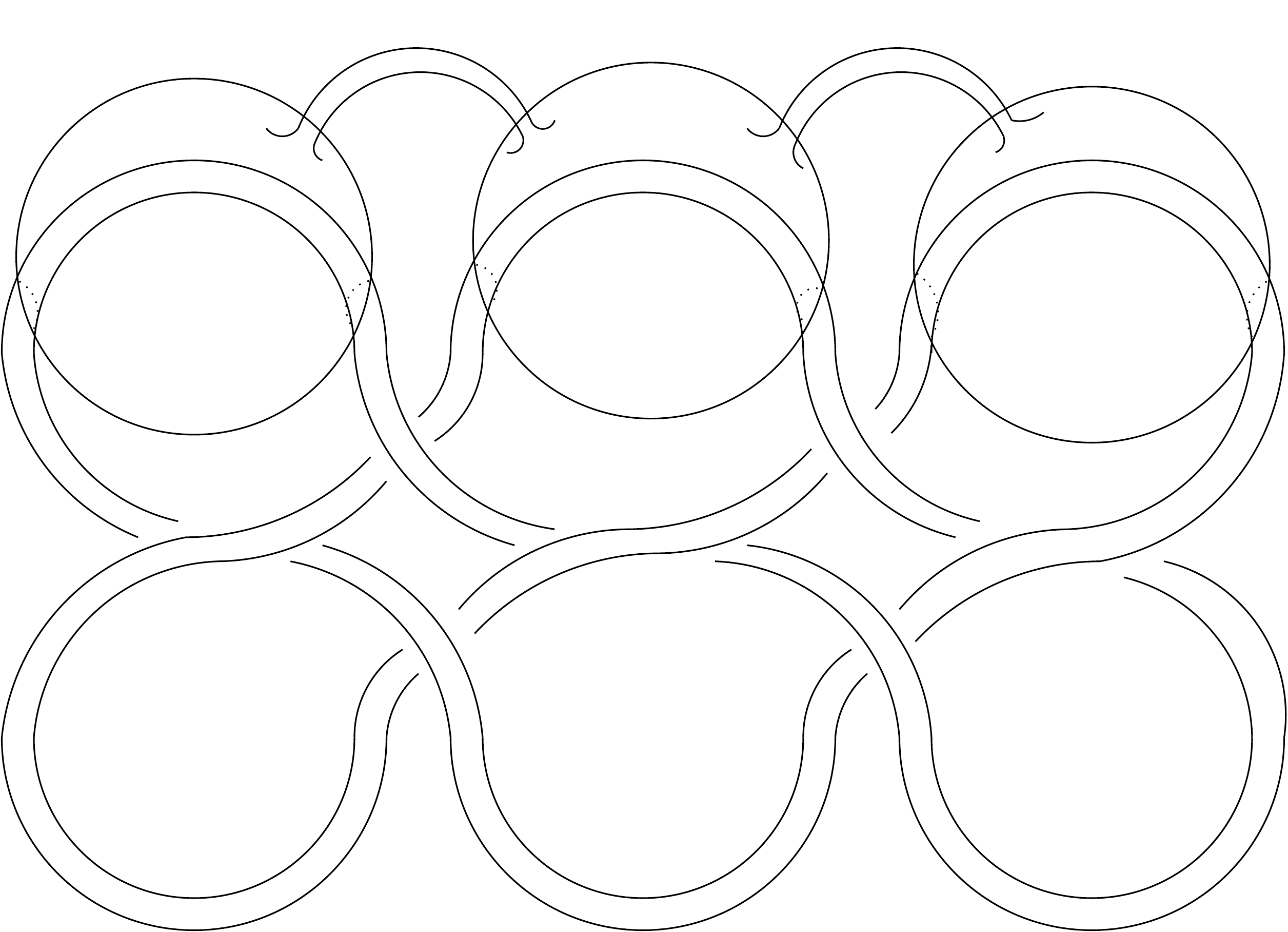}}
\caption{The result of isotoping the annuli $A_1, A_2, A_3$ into $M_2^3$}
\label{dsa11}
\end{figure}

Consider the portion of $F^n$ lying in $M_2^n$.  See \fullref{dsa11}.  It is a punctured sphere.  Moreover, it is isotopic to a
punctured sphere that consists of a level disk with $2n$ punctures and
an upper hemisphere.  See \fullref{dsa12}.  Now note that the
portion of ${\mathbb S}^3$ above a bridge sphere that coincides with
this level punctured disk and above the upper hemisphere is a 3--ball.
(Replacing the upper hemisphere of this sphere with a level disk is
equivalent to isotoping the upper hemisphere of this sphere through
infinity.  For details, see \cite[Lemma 1]{Sch1}.)  Thus the portion
of $F^n$ lying in $M_2^n$ is isotopic to a bridge sphere.  It is
hence as required in $M_2^n$.

\begin{figure}[ht!]
\centerline{\includegraphics[width=0.7\textwidth]{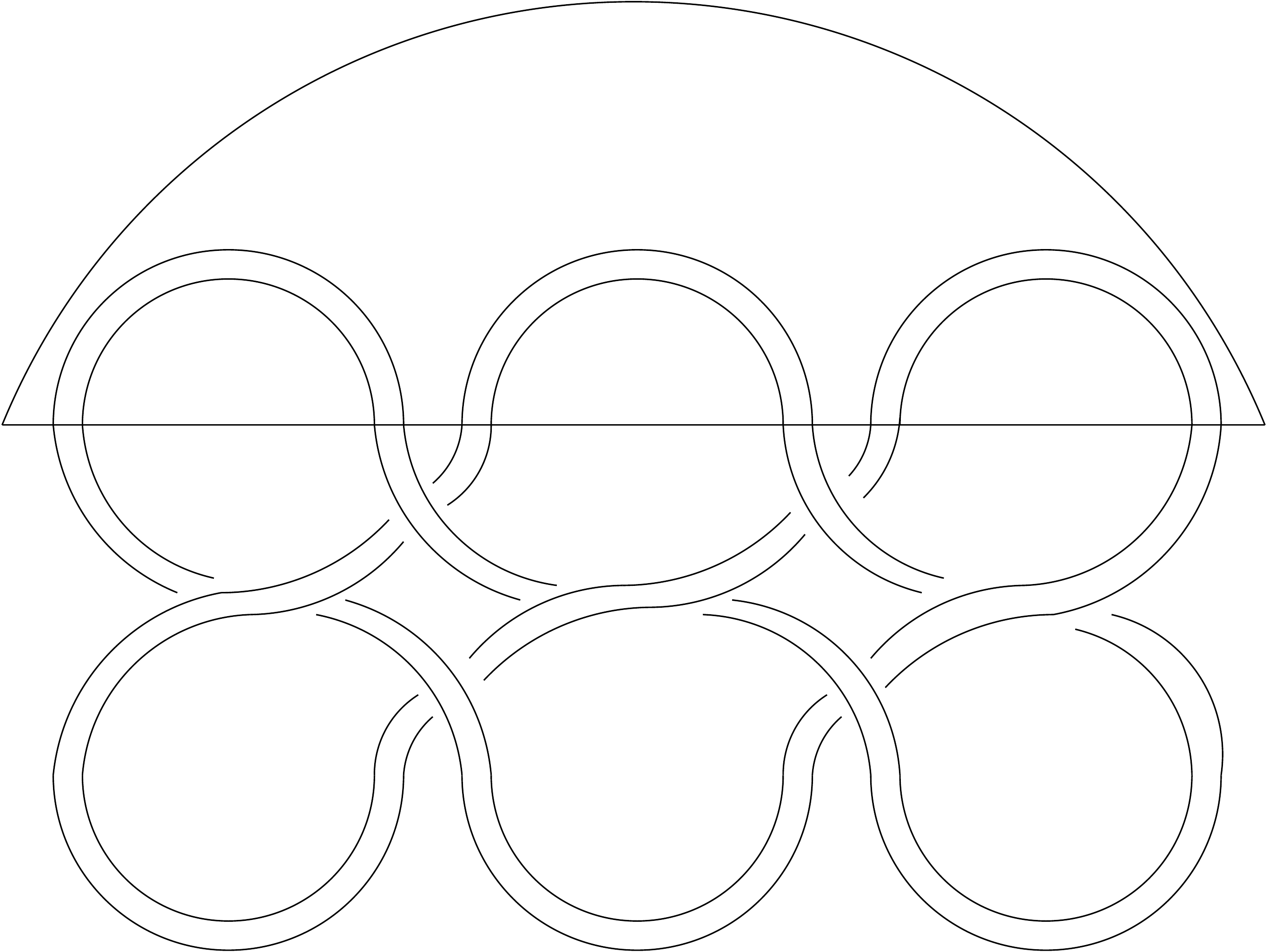}}
\caption{The punctured sphere in $M_2^3$ that is isotopic to a
bridge sphere}
\label{dsa12}
\end{figure}

It now suffices to verify that the portion of $F^n$ lying in $M_1^n$
admits the required pairs of disks.  After a small isotopy, $b_1,
\dots, b_n$ lie in the interior of $M_1^n$.  We then see that the
portion of $F^n$ lying in $M_1^n$ may be reconstructed from $n$
vertical annuli and one torus by ambient 1--surgery along arcs dual to
$b_1, \dots, b_n$.  See \fullref{dsa13}.  (Compare to \fullref{dsa10}.)

\begin{figure}[ht!]
\centerline{\includegraphics[width=0.5\textwidth]{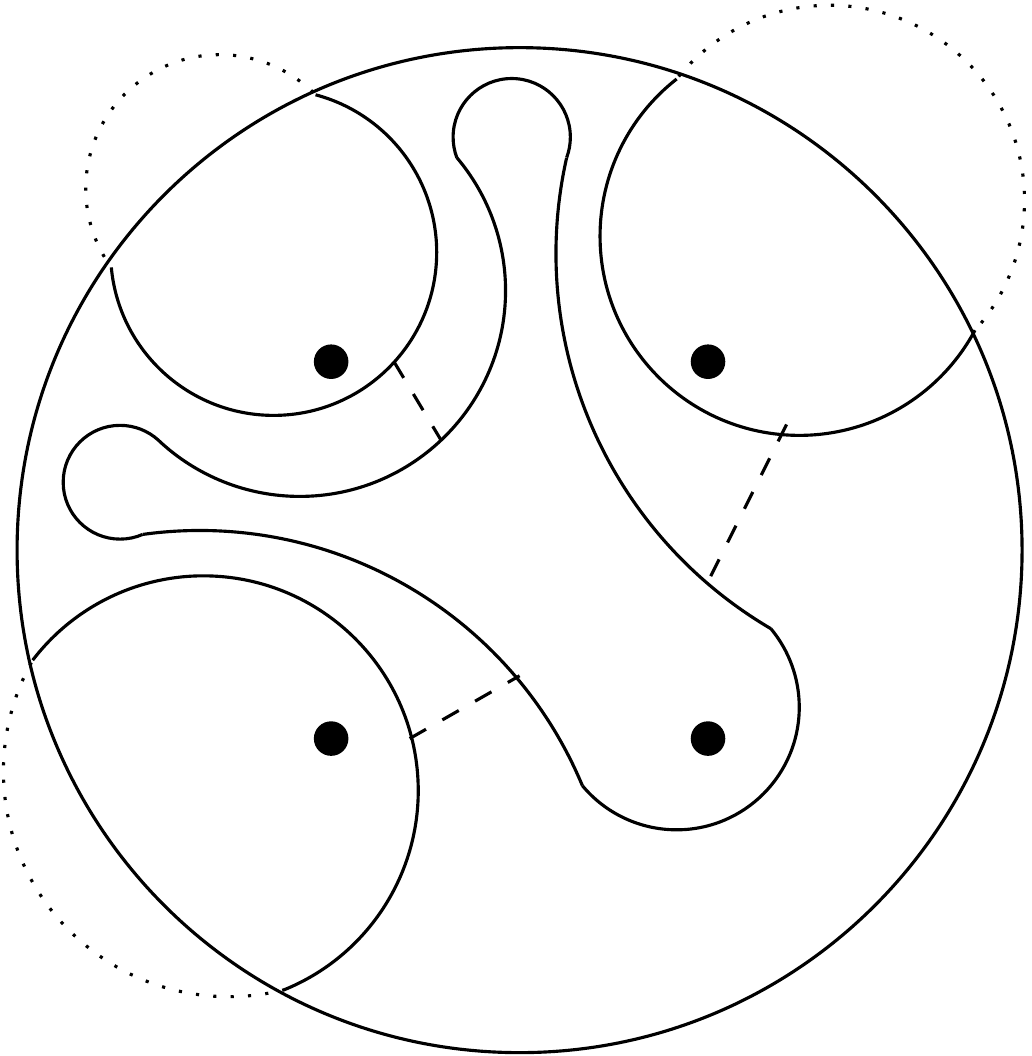}}
\caption{A dual schematic for $F^n \cap M_1^n$}
\label{dsa13}
\end{figure}

Comparing the decomposition here with $(V^n, W^n)$, we see that the
splitting surface $F^n$ is entirely contained in a collar of one of
the handlebodies $V^n$, $W^n$, say $V^n$.  Furthermore, it induces a
Heegaard splitting $(X_v^n, Y_v^n)$ of $V^n$ as follows: Take $X_v^n$
to be $X^n \cap V^n = X^n$ and take $Y_v^n$ to be the collar of
$\partial V^n$ together with $Y^n \cap V^n$.  Then $Y_v^n = (\text{collar
of } V^n) \cup (\text{solid torus}) \cup (\text{1--handles})$ and hence $X_v^n$
and $Y_v^n$ are both handlebodies.

However, the genus of $F^n$ is $2n$ and the genus of $\partial V^n$
is $n$.  It thus follows from Scharlemann and Thompson \cite[Lemma
2.7]{ST} that $(X_v^n, Y_v^n)$ and thus $(X^n, Y^n)$ is stabilized.
By applying \cite[Lemma 2.7]{ST} to locate a stabilizing pair of disks
and using one of the disks to destabilize $n$ times in succession,
we locate the $n$ pairs of stabilizing disks required.
\end{proof}

\bibliographystyle{gtart}
\bibliography{link}

\end{document}